\documentclass[a4paper,12pt,reqno]{amsart}
\usepackage{amssymb}
\usepackage{ifthen}
\nonstopmode
\numberwithin{equation}{section}
\newtheorem{thm}{Theorem}[section]
\newtheorem{cor}[thm]{Corollary}
\newtheorem{lem}[thm]{Lemma}
\newtheorem{prop}[thm]{Proposition}

\theoremstyle{remark}
\newtheorem{rem}[thm]{Remark}

\newcounter{alphabet}
\newcounter{tmp}

\newenvironment{pf}[1][]{%
 \vskip 3mm
 \noindent
 \ifthenelse{\equal{#1}{}}%
  {{\slshape Proof. }}%
  {{\slshape #1.} }%
 }%
{\qed\bigskip}

\newcounter{minutes}
\divide\time by 60
\newcounter{hours}
\multiply\time by 60
\addtocounter{minutes}{-\time}
\begin{document}
\newcommand{\A}{{\mathcal A}}
\newcommand{\B}{{\mathcal B}}
\newcommand{\T}{{\mathcal T}}
\newcommand{\M}{{\mathcal M}}
\newcommand{\F}{{\mathcal F}}
\newcommand{\Om}{{\Omega}}
\newcommand{\es}{{\mathcal S}}
\newcommand{\R}{{\mathbb R}}
\newcommand{\C}{{\mathbb C}}
\newcommand{\K}{{\mathcal K}}
\newcommand{\eK}{{\bold K}}
\newcommand{\uhp}{{\mathbb H}}
\newcommand{\Z}{{\mathbb Z}}
\newcommand{\N}{{\mathbb N}}
\newcommand{\D}{{\mathbb D}}
\newcommand{\UCV}{{\mathcal{UCV}}}
\newcommand{\kUCV}{{k\text{-}\mathcal{UCV}}}
\newcommand{\Perron}{{\mathcal P}}
\newcommand{\sphere}{{\widehat{\mathbb C}}}
\newcommand{\image}{{\operatorname{Im}\,}}
\renewcommand{\Im}{{\operatorname{Im}\,}}
\newcommand{\Aut}{{\operatorname{Aut}\,}}
\newcommand{\real}{{\operatorname{Re}\,}}
\renewcommand{\Re}{{\operatorname{Re}\,}}
\newcommand{\kernel}{{\operatorname{Ker}\,}}
\newcommand{\id}{{\operatorname{id}}}
\newcommand{\mob}{{\text{\rm M\"{o}b}}}
\newcommand{\Int}{{\operatorname{Int}\,}}
\newcommand{\Ext}{{\operatorname{Ext}\,}}
\renewcommand{\mod}{{\operatorname{mod}\,}}
\newcommand{\stab}{{\operatorname{Stab}}}
\newcommand{\SL}{{\operatorname{SL}}}
\newcommand{\PSL}{{\operatorname{PSL}}}
\newcommand{\PSU}{{\operatorname{PSU}}}
\newcommand{\tr}{{\operatorname{tr}}}
\newcommand{\diam}{{\operatorname{diam}\,}}
\newcommand{\inv}{^{-1}}
\newcommand{\area}{{\operatorname{Area}\,}}
\newcommand{\eit}{{e^{i\theta}}}
\newcommand{\eint}{{e^{in\theta}}}
\newcommand{\emint}{{e^{-in\theta}}}
\newcommand{\dist}{{\operatorname{dist}}}
\newcommand{\arctanh}{{\operatorname{arctanh}}}
\newcommand{\const}{{\operatorname{const.}}}
\newcommand{\capa}{{\operatorname{Cap}}}
\newcommand{\hdim}{{\operatorname{H-dim}}}
\newcommand{\rad}{{\operatorname{rad}}}
\newcommand{\partialb}{{\partial_{\operatorname{b}}}}
\newcommand{\CD}{{\operatorname{CD}}}
\newcommand{\hm}{{\mathcal H}}
\newcommand{\hc}{{\mathcal L}}
\newcommand{\cube}{{\mathcal Q}}
\newcommand{\Log}{{\,\operatorname{Log}}}

\newcommand{\Der}{{\frak D}}
\newcommand{\Isom}{{\operatorname{Isom}}}
\newcommand{\der}{{\mathcal D}}
\renewcommand{\r}{{\varphi}}
\newcommand{\s}{{\psi}}
\newcommand{\hol}{{\mathrm H}}

\bibliographystyle{amsplain}
\title%
[On a Multi-point Schwarz-Pick Lemma]%
{On a Multi-point Schwarz-Pick Lemma}

\author[K.~H.~Cho]{Kyung Hyun Cho}
\address{Department of Physics, Pohang University of Science and Technology,
Pohang, Kyungbuk 790-784, Korea}
\email{khsacho@postech.ac.kr}

\author[S.-A Kim]{Seong-A Kim}
\address{Department of Mathematics Education, Dongguk University \\
Gyeongju, Kyungbuk 780-714, Korea}
\email{sakim@dongguk.ac.kr}  

\author[T.~Sugawa]{Toshiyuki Sugawa}
\address{Graduate School of Information Sciences, Tohoku University \\
Aoba-ku, Sendai 980-8579, Japan}
\email{sugawa@math.is.tohoku.ac.jp}

%\date{}

\begin{abstract}
We consider the multi-point Schwarz-Pick lemma and its associate
functions due to Beardon-Minda and Baribeau-Rivard-Wegert.
Basic properties of the associate functions are summarized.
Then we observe that special cases of the multi-point Schwarz-Pick lemma give
Schur's continued fraction algorithm and several inequalities
for bounded analytic functions on the unit disk.
\end{abstract}

\thanks{
The third author was supported in part by JSPS Grant-in-Aid for Scientific
Research (B), 17340039 and for Exploratory Research, 19654027.
}

%\begin{center}
%{\tiny \texttt{FILE:~\jobname .tex,
%        printed: \number\year-\number\month-\number\day, 
%        \thehours.\ifnum\theminutes<10{0}\fi\theminutes}
%}
%\end{center}

\maketitle

\section{Introduction and preliminaries}
Many ways of applying the Schwarz lemma reveal deep properties of
holomorphic mappings $f:\D\to\D,$ where $\D$ will denote the unit disk 
$\{z\in\C:|z|<1\}$ throughout the present paper.
For instance, the refined forms of the Schwarz Lemma due to Dieudonn\'e and
Rogosinski are explained in detail in \cite{Duren:univ}.
More recently, a number of sharpened forms of the Schwarz or
Schwarz-Pick Lemma have been obtained 
(see \cite{BRW09}-\cite{BM08} and \cite{Mercer99}).
Among others, Beardon and Minda \cite{BM04} presented 
an extension of the Schwarz-Pick Lemma which involves three points 
and yields known variations of the Schwarz-Pick Lemma in a unified way.
Later on, Baribeau, Rivard and Wegert \cite{BRW09} generalized it to $n$ points
and applied it to Nevanlinna-Pick interpolation problem.

In this paper, we discuss the multi-point Schwarz-Pick Lemma by defining
a set of holomorphic functions on $\D$ associated 
with a sequence of given points in  $\D.$
We observe how our results are related with the Schur algorithm and
show that they turn to coefficient estimates for a bounded analytic function on  $\D$ 
and there is a correlation between the coefficient estimates.
Moreover, we obtain some applications of the results.
We now start by recalling the Schwarz-Pick Lemma.

\begin{lem}\label{lem:sp}
Let $f:\D\to\D$ be holomorphic and fix $z_0\in\D.$
For any point $z \in \D$, the inequality
\begin{equation}\label{eq:sp}
\left|\frac{f(z)-f(z_0)}{1-\overline{f(z_0)}f(z)}\right| 
\le \left|\frac{z-z_0}{1-\overline{z_0}z}\right|
\end{equation}
holds if $z\ne z_0$ and
\begin{equation}\label{eq:spd}
\frac{|f'(z)|}{1-|f(z)|^2} \le \frac1{1-|z|^2}
\end{equation}
if $z=z_0.$
Equality holds for a point $z$ precisely when 
$f$ is a conformal automorphism of the unit disk $\D.$
\end{lem}

We denote by $\Aut(\D)$ the group of conformal automorphisms of $\D.$
Note that $f:\D\to\D$ is in $\Aut(\D)$ if and only if 
$f(z)=(\alpha z+\beta)/(\bar\beta z+\bar\alpha)$ for complex constants
$\alpha$ and $\beta$ with $|\alpha|^2-|\beta|^2=1.$
Furthermore, $\hol(\D)$ will denote the set of holomorphic functions $f$
on $\D$ with $|f|\le 1.$
By the maximum principle, $f\in\hol(\D)$ is a (unimodular) constant
if $f$ assumes a value in the boundary  $\partial\D$ of  $\D.$
In other words, $|f|<1$ in $\D$ for $f\in\hol(\D)$ unless $f$ is a constant.

For $ z, w\in\D,$ let $[z,\,w]$ be defined by
\begin{equation}\label{eq:dist}
[z,\,w] = \frac{z-w}{1-\overline{w}z}.
\end{equation}

Its modulus $|[z,\,w]|$ is called the pseudo-hyperbolic distance between
$z$ and $w$ in $\D$ \cite{BM04}.
It is convenient to memorize the fact that $w=[z,z_0]$ if and only if
$z=[w,-z_0]$ for three points $z_0, z, w\in\D.$
We extend the definition of $[z,\,w]$ by letting
$[z,\,z]=0$ for $z\in\partial\D$ so that $[f(z),\,f(w)]$ is defined
whenever $f\in\hol(\D)$ and $z,w\in\D.$
The inequality \eqref{eq:sp} is now same as
\begin{equation}\label{eq:sp-psdist}
|[f(z),\,f(z_0)]| \le |[z,\,z_0]|.
\end{equation}
The geometrical meaning of the Schwarz-Pick Lemma is that $f$ is 
distance-decreasing with respect to the hyperbolic
metric $\rho(z)|dz|=2|dz|/(1-|z|^2)$ of the unit disk.
We denote by $d(z,w)$ the hyperbolic distance induced by $\rho;$
in other words,
$$
d(z,w)=\log\frac{1+|[z,\,w]|}{1-|[z,\,w]|}.
$$
The inequality \eqref{eq:sp} is equivalent to $d(f(z),f(z_0))\le d(z,z_0)$
for a holomorphic map $f:\D\to\D$ and $z, z_0\in\D.$

\bigskip

Let us briefly recall the main idea of Beardon and Minda \cite{BM04}.
For this purpose, we introduce an operation for functions as follows.
Let $f\in\hol(\D)$ and $z_0\in\D.$
We define a holomorphic function $\Delta_{z_0}f$ on $\D$ by
\begin{equation}\label{eq:hdquo}
\Delta_{z_0}f(z)=
\begin{cases}
\dfrac{[f(z),\,f(z_0)]}{[z,\,z_0]} &\quad \text{for}~z\ne z_0, \\
\null & \\
\dfrac{(1-|z_0|^2)f'(z_0)}{1-|f(z_0)|^2} &\quad\text{for}~z=z_0.
\end{cases}
\end{equation}
The symbol $\Delta_{z_0}f$ is adopted in \cite{BRW09}.
When it is convenient to regard $\Delta_{z_0}f(z)$ as a function of
the two variables $z$ and $z_0,$ we also write
$\Delta_{z_0}f(z)=f_1(z; z_0).$
In \cite{BM04}, this quantity is called the 
`hyperbolic difference quotient' of $f,$ and
the above notation is
somewhat different from that of \cite{BM04} for the purpose
of introducing hyperbolic difference quotients of higher order.
By the form of the definition, we have naturally the chain rule
$$
\Delta_{z_0}(f\circ g)
=(\Delta_{g(z_0)}f)\circ g\cdot \Delta_{z_0}g
$$
for holomorphic maps $f, g:\D\to\D$ and $z_0\in\D$ (cf.~\cite{BM04}).
Since $\Delta_{z_0}T(z)=T'(z_0)/|T'(z_0)|$ for $T\in\Aut(\D),$
the following invariance property can easily be deduced.

\begin{lem}\label{lem:inv}
Let $f:\D\to\D$ be a holomorphic map.
For conformal automorphisms $S$ and $T$ of $\D$,
\begin{equation*}%\label{eq:inv1}
\Delta_{z_0}(S \circ f \circ T)(z) =
\frac{S'(f(T(z_0)))}{|S'(f(T(z_0)))|}
\cdot\frac{T'(z_0)}{|T'(z_0)|}\cdot \Delta_{T(z_0)}f(T(z))
%(S \circ f \circ T)_1 (z;z_0) =
%f_1(T(z);T(z_0))\frac{|S'(f(T(z_0)))|}{S'(f(T(z_0)))}
%\cdot\frac{|T'(z_0)|}{T'(z_0)} 
\end{equation*}
for $z, z_0\in\D.$
\end{lem}

In particular, $|(S \circ f \circ T)_1 (z;z_0)| =|f_1(T(z);T(z_0))|$
%and $|D_1(S \circ f \circ T) (z)| =|D_1f(T(z))|$
(cf.~\cite[Lemma 2.3]{BM04}).

In terms of the hyperbolic difference quotient,
the Schwarz-Pick Lemma is now rephrased as follows.

\begin{lem}\label{lem:f1}
Let $f\in\hol(\D).$
Then, for any pair of points $z, z_0 \in \D$,
\begin{equation}\label{eq:sp-hdquo}
|\Delta_{z_0}f(z)|\equiv|f_1(z;z_0)| \le 1.
\end{equation}
Here, equality holds for a pair of points 
precisely when $f\in\Aut(\D).$
\end{lem}

Note that $\Delta_{z_0}f$ is a unimodular constant for any $z_0\in\D$ when 
$f\in\Aut(\D)$ and that $|\Delta_{z_0}f|<1$ on $\D$ when 
$f\in\hol(\D)\setminus\Aut(\D).$

It is crucial to note that, by the Schwarz-Pick Lemma,
the function $\Delta_{z_0}f$ again belongs to $\hol(\D)$ for $f\in\hol(\D)$
and $z_0\in\D;$
in other words, $\Delta_{z_0}$ is an operator on $\hol(\D)$ into itself.
This observation leads to the following definition (cf.~\cite{BRW09}):
Let $f\in\hol(\D).$
For a given (finite or infinite) sequence of points $z_j~(j=0,1,\dots)$ 
in $\D,$ define $f_j(z; z_{j-1}, \dots, z_0)~ (j=0,1,\dots)$ by
$$
f_{j}(z; z_{j-1},\dots,z_0)=(\Delta_{z_{j-1}}\circ\cdots\circ\Delta_{z_0})f(z).
$$
Here, we understand that $f_0(z; -)=f(z)$ for $j=0.$
Note that this notation is consistent with the former definition
of $f_1(z; z_0).$

For brevity, we also write $f_j(z)=f_j(z; z_{j-1},\dots,z_0)$
and $\gamma_j=f_j(z_j)$ for $j=0,1,2,\dots.$
%The sequence $\gamma_j$ is called the Schur parameters.
%We have then two possibilities for the Schur parameters.
We have then two possibilities:

\noindent
(i) $|\gamma_j| < 1$ for each $j.$
Then $|f_j|<1$ for each $j.$
If $f_j$ is constant for some $j,$ then $f_k = \gamma_k = 0$ for $k>j.$

\noindent
(ii) There exists an integer $n$ such that 
$|\gamma_0| < 1, |\gamma_1| < 1, \cdots, 
|\gamma_{n-1}| < 1, |\gamma_n| = 1$. 
Then, $f_n=\gamma_n$ and $f$ turns out to be a Blaschke product of degree $n.$
Beardon and Minda \cite{BM04} showed that this occurs only in this case.
Here, we recall that a function $f$ is called a (finite) Blaschke product
of degree $n$ if $f(z)=e^{i\theta}\prod_{j=1}^n[z,a_j]$ for $\theta\in\R$ and
some points $a_1,\dots,a_n\in\D.$
Note that $f_j=0$ for $j>n$ in this case.

Through the above observation, for $f\in\hol(\D),$
we see that $|f_n(z)|=1$ for some $z\in\D$
if and only if $f$ is a Blaschke product of degree $n.$

%Remark: Schur algorithm (moved to \S 3)

By repeated applications of the Schwarz-Pick lemma, 
we now have the following multi-point Schwarz-Pick Lemma due to
Beardon-Minda \cite{BM04} for $j=2$ and Baribeau-Rivard-Wegert \cite{BRW09}
for general $j.$

\begin{thm}\label{thm:n-sp}
Let $f\in\hol(\D)$ and $z_0, z_1, \dots, z_j$ be a sequence of 
$j+1$ points in $\D.$
Then
\begin{equation}\label{eq:sp-n}
|f_j(z;z_{j-1},\dots,z_0)|\le 1,\quad z\in\D.
\end{equation}
Equality holds for a point $z\in\D$ if and only if $f$ is a Blaschke product
of degree $j.$
Moreover, if $f$ is not a Blaschke product of degree $\le j,$
\begin{equation}\label{eq:sp-disti}
d(f_j(z; z_{j-1}, \cdots, z_0),\,f_j(z_j; z_{j-1}, \cdots, z_0)) \le d(z,\,z_j),
\quad z\in\D.
\end{equation}
Equality holds for a point $z\ne z_j$ precisely
when $f$ is a Blaschke product of degree $j+1.$
\end{thm}

It is shown in \cite{BM04} that many known results 
in \cite{Duren:univ} and \cite{Mercer99}
can be derived based on the above theorem for $j=2;$ namely
the `three-point' Schwarz-Pick Lemma. 
In the present note, we give some consequences of $n$-point Schwarz-Pick
Lemma. To this end, we also present a couple of basic properties of
the quantities $f_j(z_j; z_{j-1},\dots, z_0)$ for $f$ and $z_0,\dots, z_j$
in the next section.
In Section 3, several interpretations and applications are given.
Indeed, we will point out relations to the Schur algorithm and Peschl's
invariant derivatives, and give several concrete refinements
of known results such as Yamashita's inequality, Dieudonn\'e's lemma.
For Dieudonn\'e's lemma \cite{Duren:univ}, in addition to \cite{BM04},
see also \cite{BM08}.
We would like to remark that such refinements could be given, in principle,
as much as we wish, with the expense of complication.

\section{Main Results}\label{sec:main}

We first observe analyticity of the function $f_j(z; z_{j-1}, \dots, z_0)$
for $f\in\hol(\D).$
This property guarantees existence of the limit of $f_j(z_j; z_{j-1}, \dots, z_0)$
as $z_k\to z_l$ for a pair of the variables $z_k$ and $z_l$ for instance,
and allows us to change the order of limits.
% concerning the variables $z_0, \dots, z_j$ and $z.$

\begin{prop}
Let $f\in\hol(\D).$
Then for each $j\ge0,$ the function $f_j(z; z_{j-1}, \dots, z_0)$
is complex analytic in $z\in\D$ and real analytic in $z_0, \dots, z_{j-1}\in\D.$
\end{prop}

\begin{pf}
We show the assertion by induction on $j.$
It is clear for $j=0.$
We assume that the assertion is valid up to $j.$
By definition,
\begin{align*}
&f_{j+1}(z; z_j,\dots,z_0) \\
&=\frac{f_j(z;z_{j-1},\dots,z_0)-f_j(z_j;z_{j-1},\dots,z_0)}{z-z_j}
\cdot
\frac{1-\overline{z_j}z}%
{1-\overline{f_j(z_j;z_{j-1},\dots,z_0)}f_j(z;z_{j-1},\dots,z_0)}.
\end{align*}
There is nothing to show when $z\ne z_j.$
Thus, it is enough to show analyticity at every point of the form
$(z,z_j,z_{j-1},\dots,z_0)=(a_j,a_j,a_{j-1},\dots,a_0).$

The second factor of the right-hand side in the above formula
is clearly analytic in the sense of the assertion.
Analyticity of the first factor follows from the next lemma by
the interpretation $z_k=t_{2k+1}+it_{2k+2}$ for $k=0,1,\dots, j-1$
and $w=z_j.$
\end{pf}

\begin{lem}
Suppose that a continuous function $F(z,t_1,\dots,t_n)$ is
complex analytic in the complex variable $z$ and real analytic
in the real variables $t_1,\dots,t_n.$
Then the difference quotient
$$
\frac{F(w,t_1,\dots,t_n)-F(z,t_1,\dots,t_n)}{w-z}
$$
is complex analytic in $z, w$ and real analytic in $t_1,\dots,t_n.$
\end{lem}

\begin{pf}
For simplicity, we prove only in the case when $n=1.$
It is enough to see that $(F(w,t)-F(z,t))/(w-z)$ is complex analytic
in $|z|<r/2, |w|<r/2$ and real analytic in $|t|<\delta$
for small enough $r>0$ and $\delta>0.$
We may assume that $F$ is expanded in the form
$$
F(z,t)=\sum_{j=0}^\infty A_j(z)t^j,\quad |z|<2r, |t|<2\delta
$$
for some constants $r>0$ and $\delta>0.$

By convergence of the above series, 
there exists a constant $M>0$ such that
$$
|A_j(z)|\le M\delta^{-j},\quad |z|\le r,~ j=0,1,2,\dots.
$$
Since $F(z,t)$ is complex analytic in $z,$ Cauchy's integral formula
yields the expression
\begin{align*}
\frac{F(w,t)-F(z,t)}{w-z}
&=\frac1{2\pi i}\int_{|\zeta|=r}\frac{F(\zeta,t)}{(\zeta-w)(\zeta-z)}d\zeta \\
&=\sum_{j=0}^\infty 
\frac1{2\pi i}\int_{|\zeta|=r}
\frac{A_j(\zeta)}{(\zeta-w)(\zeta-z)}d\zeta\cdot t^j \\
&\equiv\sum_{j=0}^\infty B_j(z,w)t^j
\end{align*}
for $|z|<r, |w|<r.$
Here,
$$
|B_j(z,w)|\le \frac{4M}{r\delta^j},\quad |z|<r/2,~ |w|<r/2,
$$
and thus the above series is indeed convergent in
$|t|<\delta.$
\end{pf}

The following generalization of Lemma \ref{lem:inv} will be useful
to reduce general questions to special ones.

\begin{lem}\label{lem:invg}
Let $f\in\hol(\D), S, T\in\Aut(\D)$ and $z_0,\dots, z_{j-1}, z\in\D.$
Then
$$
(S\circ f\circ T)_j(z; z_{j-1},\dots, z_0)
=\frac{S'(f(T(z_0)))}{|S'(f(T(z_0)))|}
\cdot f_j(T(z); T(z_{j-1}),\dots, T(z_0))
\cdot\prod_{k=0}^{j-1}\frac{T'(z_k)}{|T'(z_k)|}.
$$
In particular,
$$
|(S\circ f\circ T)_j(z; z_{j-1},\dots, z_0)|
=|f_j(T(z); T(z_{j-1}),\dots, T(z_0))|.
$$
\end{lem}

\begin{pf}
We can easily verify the relation $[\zeta z,\zeta w]=\zeta[z,w]$
for $z,w\in\D$ and $\zeta\in\partial\D.$
Therefore, $$
\Delta_{z_0}(\zeta f)=\zeta\Delta_{z_0}f
$$
for $f\in\hol(\D), z_0\in\D$ and $\zeta\in\partial\D.$
For brevity, we put $\omega=S'(f(T(z_0)))/|S'(f(T(z_0)))|,
\zeta_k=T'(z_k)/|T'(z_k)|,~z'=T(z)$ and $z_k'=T(z_k).$
Note that $\omega,\zeta_k\in\partial\D.$
By Lemma \ref{lem:inv} together with the above relation, we see
\begin{align*}
(S\circ f\circ T)_2(z; z_1,z_0)
&=\Delta_{z_1}(\Delta_{z_0}(S\circ f\circ T))(z) \\
&=\Delta_{z_1}(\omega\zeta_0 (\Delta_{z_0'}f)\circ T)(z) \\
&=\omega\zeta_0\Delta_{z_1}((\Delta_{z_0'}f)\circ T)(z) \\
&=\omega\zeta_0\zeta_1\Delta_{z_1'}(\Delta_{z_0'}f)(T(z)) \\
&=\omega\zeta_0\zeta_1f_2(z'; z_1', z_0').
\end{align*}
In the same way, we can show the required relation for general $j.$
\end{pf}

Let $f\in\hol(\D)$ and a point $z\in\D$ be given.
The most crude estimate for $f(z)$ is $|f(z)|\le 1.$
This, however, cannot be improved without any additional information about $f.$
If we know about the value $w_0$ of $f$ at a given point $z_0,$ then
the estimate can be improved.
For instance, when $f(0)=w_0,$ we have the better estimate \cite[p.~167]{Nehari:conf}
\begin{equation}\label{eq:nehari}
\frac{|f(0)|-|z|}{1-|z||f(0)|} \le |f(z)| \le  \frac{|f(0)|+|z|}{1+|z||f(0)|}.
\end{equation}
When more values of $f$ (and possibly its derivatives) at 
points $z_j$ for $j=0, 1, 2,  \cdots$ are specified, we may
improve the estimate more.
Indeed, we are able to show the following.
%As the reader will see, the proof is almost same as in
%Nevanlinna \cite{Nev19} (see also Theorem \ref{thm:NP} below).

\begin{thm}\label{thm:mob}
Let $a, z_0,\dots, z_n$ be given points in $\D$
and put $\tau_j=[a,z_j]$ for $j=0,1,\dots,n.$
\begin{enumerate}
\item[(i)]
Suppose that $f\in\hol(\D)$ is not a Blaschke product of
degree at most $n.$
Let $f_j(z)=f_j(z; z_{j-1},\dots, z_0),$
$\gamma_j=f_j(z_j)$ for $j=0,1,\dots,n.$
Define M\"obius transformations $A_j,~j=0,1,\dots,n,$ by
$$
A_j(x)=\frac{\tau_j x+\gamma_j}{1+\overline{\gamma_j}\tau_j x}.
$$ 
Then $f(a)\in (A_0\circ\cdots\circ A_n)(\overline\D).$
If furthermore $f$ is not a Blaschke product of degree $n+1,$
$f(a)\in (A_0\circ\cdots\circ A_n)(\D).$
\item[(ii)]
Conversely, suppose that points $\gamma_0, \gamma_1,\dots, \gamma_n
\in\D$ are given.
Let $A_j$ be as above and choose an arbitrary point $b\in
(A_0\circ\cdots\circ A_n)(\overline\D).$
Then there exists a function $f\in\hol(\D)$ with $f(a)=b$ such that
$\gamma_j=f_j(z_j;z_{j-1},\dots,z_0)$ for $j=0,1,\dots, n.$
\end{enumerate}
\end{thm}

\begin{pf}
We first show (i).
Let $w_j=f_j(a)$ for $j=0,1,\dots, n+1.$
Here, $f_{n+1}(z)$ is defined similarly.
By assumption, $|\gamma_j|<1$ and $|w_j|<1$ for $j\le n.$
Also note that $|w_{n+1}|\le 1$ and equality holds if and only if
$f$ is a Blaschke product of degree $n+1.$
Then, by definition,
$$
w_{j+1}=\Delta_{z_j}f_j(a)=\frac{[w_j,\gamma_j]}{\tau_j},
$$
and thus, 
\begin{equation}\label{eq:wj}
w_j=[\tau_jw_{j+1},-\gamma_j]=A_j(w_{j+1})
\end{equation}
for $j=0,1,\dots, n.$
Therefore, $w_0=(A_0\circ\cdots\circ A_n)(w_{n+1})$ and (i) is proved.

We next show (ii).
Set $c=(A_0\circ\cdots\circ A_n)\inv(b).$
Then, by assumption, $c\in\overline\D.$
Let $w_{n+1}=c$ and define $w_n,w_{n-1},\dots,w_0$ inductively by \eqref{eq:wj}.
Let $f_{n+1}$ be any function in $\hol(\D)$ such that
$f_{n+1}(a)=c.$
For instance, $f_{n+1}$ can be taken to be the constant function $c.$
Then, define functions $f_n, f_{n-1}, \dots, f_0$ inductively by the formula
\begin{equation}\label{eq:fj}
f_j(z)=[[z,z_j]f_{j+1}(z),-\gamma_j]
=\frac{[z,z_j]f_{j+1}(z)+\gamma_j}{1+\bar \gamma_j[z,z_j]f_{j+1}(z)}.
\end{equation}
Then $f_j(z_j)=[0,-\gamma_j]=\gamma_j$ and therefore
the relation $\Delta_{z_j}f_j=f_{j+1}$ holds.
We now set $f=f_0$ so that $f_j(z; z_{j-1},\dots,z_0)=f_j(z).$
In particular, $f_j(z_j;z_{j-1},\dots,z_0)=f_j(z_j)=\gamma_j.$
By \eqref{eq:fj}, we have
$f_j(a)=[\tau_jf_{j+1}(a),-\gamma_j]=A_j(f_{j+1}(a)).$
Hence, $f(a)=f_0(a)=(A_0\circ\cdots\circ A_n)(f_{n+1}(a))=b.$
Thus, we have shown the existence of such an $f.$
\end{pf}

In applications of the last theorem, it is convenient to note the following
elementary fact: For a M\"obius transformation $A(z)=\frac{az+b}{cz+d}$
with $|c|<|d|,$
\begin{equation}\label{eq:A}
w\in A(\overline\D) \Leftrightarrow
\left|w-\frac{a\bar c-b\bar d}{|c|^2-|d|^2}\right|\le 
\left|\frac{ad-bc}{|c|^2-|d|^2}\right|.
\end{equation}
For instance, $f(a)\in A_0 (\overline\D)$ in the theorem means the inequality
$$
\left| f(a)-\frac{(1-|\tau_0|^2)\gamma_0}{1-|\gamma_0\tau_0|^2}\right|
\le\frac{(1-|\gamma_0|^2)|\tau_0|}{1-|\gamma_0\tau_0|^2},
$$
where $\tau_0=[a,z_0]$ and $\gamma_0=f(z_0).$

As another application of the relation \eqref{eq:wj}, we obtain the
next result.

\begin{thm}\label{thm:dineq}
Let $f\in\hol(\D)$ and $z_0\in\D.$
Then the double inequality
$$
\frac{||f(z_0)|-|[z,z_0]f_1(z; z_0)||}{1-|[z,z_0]f(z_0)f_1(z;z_0)|}\le |f(z)|\le
\frac{|f(z_0)|+|[z,z_0]f_1(z;z_0)|}{1+|[z,z_0]f(z_0)f_1(z;z_0)|}
$$
holds for $z\in\D.$
Equality holds in the left-hand (right-hand) inequality if and only if
either $f(z_0)f(z)=0$ 
or else $\arg f(z)=\arg f(z_0) ~(\mod 2\pi)$ (respectively, 
$\arg f(z)=\arg f(z_0)+\pi ~(\mod 2\pi)$).
\end{thm}

\begin{pf}
We first note the elementary inequalities (cf.~\cite[p.~167]{Nehari:conf})
\begin{equation}\label{eq:ab}
\frac{||b|-|a||}{1-|ab|}
=|[|b|,|a|]|\le |[a,b]|\le [|a|,-|b|]=\frac{|a|+|b|}{1+|ab|}
\end{equation}
for $a,b\in\D.$
Here, equality holds in the left-hand (right-hand) side if and only if
either $ab=0$ or else $(a/b)>0$ (resp.~$(a/b)<0$).
We now apply the above inequality to the choice 
$a=[f(z), f(z_0)]=[z,z_0]f_1(z;z_0)=\tau_0w_1$ and $b=-f(z_0)=-\gamma_0.$
Since $[a,b]=[\tau_0w_1,-\gamma_0]=w_0=f(z)$ by \eqref{eq:wj},
we obtain the assertion.
\end{pf}

By Lemma \ref{lem:f1}, we have the following.

\begin{cor}\label{cor:ineq}
Let $f\in\hol(\D)$ and $z_0\in\D.$
Then the double inequality
$$
\max\left\{\frac{|f(z_0)|-|[z,z_0]|}{1-|[z,z_0]f(z_0)|},0\right\}
\le |f(z)|\le
\frac{|f(z_0)|+|[z,z_0]|}{1+|[z,z_0]f(z_0)|}
$$
holds for $z\in\D.$
When $z\ne z_0,$
equality holds in the right-hand side only if $f\in\Aut(\D).$
\end{cor}

%\begin{pf}
%By \eqref{eq:ab} we see also that
%\begin{equation*}\label{eq:inequality}
%\max\{[|b|,r],0\}\le |[a,b]|\le [r,-|b|],\quad |a|\le r
%\end{equation*}
%for $b\in\D$ and $0\le r<1.$
%\end{pf}

Note that the corollary reduces to \eqref{eq:nehari} when $z_0=0.$
Thus, Theorem \ref{thm:dineq} improves the inequality \eqref{eq:nehari}.

Theorem \ref{thm:mob} gives precise information about the location
of the value $f(z)$ but it might not be easy to use.
We can extract more rough but convenient estimates for $|f(z)|$ as follows.

\begin{thm}\label{thm:n-sp-ap1}
Let $a, z_0,\dots, z_n$ be given points in $\D.$
Suppose that $f\in\hol(\D)$ is not a Blaschke product of
degree at most $n.$
Put $f_j(z)=f_j(z; z_{j-1},\dots, z_0),$
$\gamma_j=f_j(z_j),~ \tau_j=[a,z_j]$ for $j=0,1,\dots,n.$
Then the chain of inequalities
\begin{equation}\label{eq:n-sp-ap1}
|f(a)|\le (T_0\circ\dots\circ T_n)(1)%\le (T_0\circ\dots\circ T_{n-1})(1)
\le \dots\le (T_0\circ T_1)(1)\le T_0(1)
\end{equation} 
hold, where $T_j$ are the functions defined by
$$
T_j(x)=\frac{|\tau_j|x+|\gamma_j|}{1+|\tau_j\gamma_j|x}.
$$ 
\end{thm}

\begin{pf} 
Let $w_j=f_j(z)$ for $j=0,1,\dots, n+1$ as before.
Note first that $T_j(x)$ is non-decreasing in $0\le x\le 1,$
that $T_j(1)\le 1,$ and that $|w_{j}|\le 1.$
Therefore, the inequalities
$$
(T_0\circ\dots\circ T_n)(1)%\le (T_0\circ\dots\circ T_{n-1})(1)
\le \dots\le (T_0\circ T_1)(1)\le T_0(1)
$$
clearly hold.
Therefore, it is enough to show the inequality
$|f(z)|\le (T_0\circ\dots\circ T_n)(1).$

By the proof of Theorem \ref{thm:mob}, we have
$f(a)=w_0=(A_0\circ\cdots\circ A_n)(w_{n+1}).$
Note that \eqref{eq:ab} implies $|A_j(w)|\le T_j(|w|)\le 1$ 
for $w\in\overline\D.$
Therefore, we have
$$
|f(a)|\le T_0(|(A_1\circ\cdots A_n)(w_{n+1})|)
\le \cdots \le (T_0\circ\cdots\circ T_n)(|w_{n+1}|)
\le (T_0\circ\cdots\circ T_n)(1),
$$
as required.
\end{pf}

The bound $T_0(1)$ in the last theorem 
is the same as in Corollary \ref{cor:ineq}.
The inequality for the next term $T_0(T_1(1))$ takes the form
$$
|f(z)|\le
\frac{|f(z_0)|+\big|\frac{z-z_0}{1-\bar z_0z}\big|
\left|\frac{f(z_1)-f(z_0)}{1-\overline{f(z_0)}f(z_1)}\right|
+\big|\frac{z-z_1}{1-\bar z_1z}\big|
\left(|f(z_0)|\left|\frac{f(z_1)-f(z_0)}{1-\overline{f(z_0)}f(z_1)}\right|
+\big|\frac{z-z_0}{1-\bar z_0z}\big|\right)}%
{1+\big|\frac{z-z_0}{1-\bar z_0z}\big|
|f(z_0)|\left|\frac{f(z_1)-f(z_0)}{1-\overline{f(z_0)}f(z_1)}\right|
+\big|\frac{z-z_1}{1-\bar z_1z}\big|
\left(\left|\frac{f(z_1)-f(z_0)}{1-\overline{f(z_0)}f(z_1)}\right|
+|f(z_0)|\big|\frac{z-z_0}{1-\bar z_0z}\big|\right)}.
$$
Since $T_0(T_1(T_2(1)))$ is too complicated to write down,
we restrict ourselves to the simple case
when $z_0=z_1=\dots=z_n=0$ so that $\tau_j=z$ for all $j.$
For brevity, we write $c_j=|\gamma_j|.$
Then the first three inequalities in Theorem \ref{thm:n-sp-ap1} can be expressed
by
\begin{align*}
|f(z)|
&\le
\frac{c_0+(c_1+c_0c_2+c_0c_1c_2)|z|+(c_0c_1+c_2+c_1c_2)|z|^2+|z|^3}%
{c_0|z|^3+(c_1+c_0c_2+c_0c_1c_2)|z|^2+(c_0c_1+c_2+c_1c_2)|z|+1} \\
&\le
\frac{c_0+(c_1+c_0c_1)|z|+|z|^2}{c_0|z|^2+(c_1+c_0c_1)|z|+1} \\
&\le
\frac{c_0+|z|}{1+c_0|z|}.
\end{align*}
Yamashita showed an inequality equivalent to the second one 
in \cite[p.~313]{Yam94P} and used it effectively
to prove uniqueness of extremal functions in the norm estimates
of starlike and convex functions of order $\alpha$ in \cite{Yam98}.
The above refinements could be used to improve the norm estimates.

As we saw before, the Schwarz-Pick lemma means the inequality
$d(f(z),f(w))\le d(z,w)$ for a holomorphic map $f:\D\to\D.$
This inequality can be refined by using the above argument.

\begin{thm}
Let $z,z_0,\dots, z_n\in\D$ and $f\in\hol(\D).$
Suppose that $f$ is not a Blaschke product of degree at most $n.$
Let $R_0(x)=(1+|\tau_0|x)/(1-|\tau_0|x)$ and 
$T_j(x)=(|\tau_j|x+|\gamma_j|)/(1+|\gamma_j\tau_j|x),$
where $\tau_j=[z,z_j]$ and $\gamma_j=f_j(z_j;z_{j-1},\dots,z_0).$
Furthermore set $R_n= R_0 \circ T_1 \circ T_2 \circ\dots \circ T_n$ for $n\ge1.$
Then,
\begin{equation}\label{eq:n-sp-ap2}
\exp(d(f(z),f(z_0))) \le R_n(1) \le R_{n-1}(1) \le \dots \le  
R_{1}(1) \le R_0(1)=\exp(d(z,z_0)).
\end{equation}
\end{thm}

\begin{pf}
Define a M\"obius transformation $S$ by $S(x)=(1+x)/(1-x).$
Then we obtain $\exp(d(f(z), f(z_0)))=S(|[w_0,\gamma_0]|)$ and 
$(S\inv\circ R_0)(x)=|\tau_0|x,$
where $w_0=f(z).$
Thus we see that \eqref{eq:n-sp-ap2} is equivalent to
$$
|\Delta_{z_0}f(z)|=\left|\frac{[w_0,\gamma_0]}{\tau_0}\right|
\le (T_1\circ\cdots\circ T_n)(1)\le\cdots\le T_1(1)\le 1,
$$
which can be obtained by applying Theorem \ref{thm:n-sp-ap1} to
the function $\Delta_{z_0}f$ and the points $z, z_1,\dots, z_n.$
\end{pf}

We consider the case when $z_0=z_1=z_2=\cdots$ 
and present explicit forms of $R_1(1)$ and $R_2(1).$
Put $t=|[z,z_0]|$ and $c_j=|f_j(z_0;z_0, \cdots, z_0)|.$
By a simple computation, we have 
$$
R_1(1)=\frac{1+2tc_1+t^2}{1-t^2}.
$$
This was first obtained in \cite{BC92}.
The improvement of this bound in the next order is
$$
R_2(1)=\frac{1+t(c_1+c_2+c_1c_2)+t^2(c_1+c_2)+t^3}%
{1+t(c_2-c_1+c_1c_2)+t^2(c_1-c_2)-t^3}.$$
Note that this is made possible by introducing 
the second order derivative of $f(z)$
through the term $c_2=|f_2(z_0;z_0,z_0)|.$

\section{Interpretations of the results and some applications}

The most immediate and potentially important application of the
multi-point Schwarz-Pick lemma is perhaps to the Nevanlinna-Pick
interpolations as was developed by Baribeau, Rivard and Wegert \cite{BRW09}.
Let us recall the Nevanlinna-Pick interpolation problem.
Let $z_0, z_1, \dots, z_n$ and $w_0, w_1, \dots, w_n$ be given points in the unit
disk $\D.$
Here, for simplicity, we assume that $z_0, \dots, z_n$ are distinct points.
%Here, we allow a pair of the points to be coincident, unlike in \cite{BRW09}.
The Nevanlinna-Pick interpolation problem asks existence of a function
$f\in\hol(\D)$ such that 
\begin{equation}\label{eq:NP}
f(z_j)=w_j \quad\text{for}\quad j=0,1,\dots, n.
\end{equation}
%(If several $z_j$'s take the same value (say $a$), 
%then the corresponding $w_j$'s must be the common value (say, $b,$) 
%and the above condition is interpreted that
%the equation $f(z)=b$ have a root $a$ with multiplicity at least
%as many as $z_j$'s.)
The solvability of the Nevanlinna-Pick interpolation problem
is characterized as positive semi-definiteness of the Hermitian form
$$
Q(t_1,\dots, t_n)=\sum_{h,k=1}^n\frac{1-w_h\bar w_k}{1-z_h\bar z_k}t_h\bar t_k
$$
(see for instance \cite[\S 1.2]{Ahlfors:conf}).

We notice that the parameters
$\gamma_0, \gamma_1, \dots, \gamma_n$ are determined only by the data
$z_0, \dots, z_n$ and $w_0, \dots, w_n$ when $f$ is a solution to
the problem \eqref{eq:NP}.
By Theorems \ref{thm:mob} and \ref{thm:n-sp-ap1}, we have the following result.

\begin{thm}\label{thm:NP}
Let $z_0, z_1, \dots, z_n$ and $w_0, w_1, \dots, w_n$ be given points in
the unit disk $\D$ with $z_j\ne z_k~(j\ne k)$
and suppose that an analytic function $f:\D\to\overline{\D}$
satisfies $f(z_j)=w_j~(j=0,1,\dots,n).$
Let $\gamma_j=f_j(z_j; z_{j-1},\dots, z_0)$ and $\tau_j=[z,z_j]$
for a fixed $z\in\D.$
Then
$$
f(z)\in (A_0\circ\cdots\circ A_n)(\overline\D)
\quad\text{and}\quad
|f(z)|\le (T_0\circ\dots\circ T_n)(1),
$$
where
$$
A_j(x)=\frac{\tau_j x+\gamma_j}{1+\overline{\gamma_j}\tau_j x},
\quad
T_j(x)=\frac{|\tau_j|x+|\gamma_j|}{1+|\tau_j\gamma_j|x}.
$$ 
\end{thm}

\begin{rem}
By the second part of Theorem \ref{thm:mob}, the set
$(A_0\circ\cdots\circ A_n)(\overline\D)$ is (so-called)
the variability region of $f(z)$ for a given $z$
concerning the solutions to the Nevanlinna-Pick interpolation problem 
in the theorem.
Since the interpolation problem does not depend on the order of the data,
this set remains unchanged if we change
the order of the interpolation data $(z_0,w_0), (z_1,w_1), \dots, (z_n,w_n).$
\end{rem}

In the previous section, we often considered the case when
$z_0=z_1=\cdots.$
This case is closely connected with the Schur algorithm and
Peschl's invariant derivatives as we now see.
Peschl's invariant derivatives
$D_nf(z),~ n=1, 2, 3, \dots,$ (with respect to the hyperbolic metric)
are defined by the series expansion
for $f\in\hol(\D)$ \cite{Peschl55} (see also \cite{KS07diff} and \cite{Schip07}):
\begin{equation*}\label{eq:invd}
[f([z,-z_0]),f(z_0)]=
\frac{f(\frac{z+z_0}{1+\bar z_0z})-f(z_0)}%
{1-\overline{f(z_0)}f(\frac{z+z_0}{1+\bar z_0z})} 
= \sum_{n=1}^{\infty} \frac{D_nf(z_0)}{n!} z^n,
%\left(\frac{z-z_0}{1-\overline{z_0}z}\right)^n,
\quad z,z_0\in\D.
\end{equation*}
Explicit forms of $D_nf(z),~n=1,2,3,$ are given by
\begin{equation*}\label{eq:d1f}
D_1f(z)=\frac{(1-|z|^2)f'(z)}{1- |f(z)|^2}, 
\end{equation*}
\begin{equation*}\label{eq:d2f}
D_2f(z)=
\frac{(1-|z|^2)^2}{1- |f(z)|^2}\left[
f''(z)-\frac{2\bar zf'(z)}{1- |z|^2}
+\frac{2\overline{f(z)}f'(z)^2}{1- |f(z)|^2}
\right],
%\frac{(1-|z|^2)^2f''(z)}{1- |f(z)|^2}
%-\frac{2\bar z(1-|z|^2)f'(z)}{1- |f(z)|^2}
%+\frac{2(1-|z|^2)^2\overline{f(z)}f'(z)^2}{(1- |f(z)|^2)^2},
\end{equation*}
and
\begin{align*}
D_3f(z)
=\frac{(1-|z|^2)^3}{1- |f(z)|^2}&\left[
f'''(z)
-\frac{6\bar zf''(z)}{1- |z|^2}
+\frac{6\overline{f(z)}f'(z)f''(z)}{1- |f(z)|^2}\right. \\
&\left. +\frac{6\bar{z}^2f'(z)}{(1-|z|^2)^2}
-\frac{12\bar z\overline{f(z)}f'(z)^2}{(1-|z|^2)(1-|f(z)|^2)}
+\frac{6\overline{f(z)}^2f'(z)^3}{(1- |f(z)|^2)^2}\right].
\end{align*}
%=\frac{(1-|z|^2)^3f'''(z)}{1- |f(z)|^2}
%+\frac{6(1-|z|^2)^3\overline{f(z)}f'(z)f''(z)}
%{(1- |f(z)|^2)^2}
%-\frac{6\bar z(1-|z|^2)^2f''(z)}{1- |f(z)|^2}
%$$
%\begin{equation*}\label{eq:d3f}
%+\frac{6\bar{z}^2(1-|z|^2)f'(z)}{1-|f(z)|^2}
%-\frac{12\bar z(1-|z|^2)^2\overline{f(z)}f'(z)^2}{(1-|f(z)|^2)^2}
%+\frac{6(1-|z|^2)^3\overline{f(z)}^2f'(z)^3}{(1- |f(z)|^2)^3}.
%\end{equation*}

Let us now recall the {\it Schur algorithm}
%(or Schur's continued fraction algorithm) 
\cite{Schur17} (see also \cite{Wall:fraction}).
Let $f\in\hol(\D).$
Define functions $f_0, f_1, f_2, \dots$ in $\hol(\D)$ inductively by
$f_0=f$ and
$$
f_{j+1}(z) = \frac{1}{z} \cdot
\frac{f_j(z)-\gamma_j}{1-\overline{\gamma_j}f_j(z)}
=\frac{[f_j(z),\gamma_j]}{[z,0]},
$$
where $\gamma_j=f_j(0).$
%which is equivalent to
%$$
%\quad f_{j}(z) = \frac{\gamma_j+zf_{j+1}(z)}{1+\overline{\gamma_j}zf_{j+1}(z)}.
%$$
The sequence $\{\gamma_j\}_{j=0}^\infty$ is called
the {\it Schur parameter} of $f.$
By construction, $f_j(z)=f_j(z; 0,\dots, 0)$ for $j=0,1,\dots.$
Recall that either $|\gamma_j|<1$ for all $j$ or else
$|\gamma_0|<1,\dots,|\gamma_{n-1}|<1, |\gamma_n|=1, \gamma_{n+1}=\cdots=0$
for some $n\ge0.$
The latter case happens precisely when $f$ is a Blaschke product of degree $n.$

We note that $D_1f(z)$
%$$D_1f(z)=\frac{(1-|z|^2)f'(z)}{1- |f(z)|^2}$$
is known as the hyperbolic derivative of $f.$
We can easily see that $f_1(z; z)=D_1f(z).$
The Schwarz-Pick lemma now implies $|D_1f(z)|\le 1.$
What is the relation between $f_n(z; z, \dots, z)$ and $D_nf(z)$?
The next result answers to it.

\begin{prop}\label{prop:rel}
Let $f\in\hol(\D)$ and $z_0\in\D.$ 
Define $g\in\hol(\D)$
by $g(z)=[f([z,-z_0]),f(z_0)].$
Then $g^{(n)}(0)=D_nf(z_0)$
and $f_n(z_0; z_0,\dots, z_0)=\gamma_n$ for $n=1,2,\dots,$ where
$\{\gamma_n\}$ is the Schur parameter of $g.$
\end{prop}

\begin{pf}
The relations $g^{(n)}(0)=D_nf(z_0)$ immediately follow from the definition
of $D_nf.$
Define $S, T\in\Aut(\D)$ by $S(w)=(w-w_0)/(1-\bar w_0w)$
and $T(z)=(z+z_0)/(1+\bar z_0z),$ where $w_0=f(z_0),$
so that $g=S\circ f\circ T.$
Note that $S'(w_0)=1/(1-|w_0|^2)>0$ and $T'(0)=1-|z_0|^2>0.$
Then, by Lemma \ref{lem:invg}, we have
$$
\gamma_n=g_n(0;0,\dots,0)=f_n(z_0;z_0,\dots,z_0).
$$
\end{pf}

When we express $g$ by the series expansion
$g(z)=\sum_{n=1}^\infty a_n z^n,$ the first several $\gamma_j$'s
are given by
\begin{align*}
\gamma_1&= a_1, \\
\gamma_2&= \frac{a_2}{1-|a_1|^2}, \\
\gamma_3&= \frac{a_3(1-|a_1|^2)+\bar a_1a_2^2}{(1-|a_1|^2)^2-|a_2|^2}, \\
\gamma_4&= \frac{a_4[(1-|a_1|^2)^2-|a_2|^2]
+2\bar a_1a_2a_3(1-|a_1|^2)+\bar a_1^2a_2^3+\bar a_2a_3^2}%
{(1-|a_1|^2)^3-(1-|a_1|^2)(|a_3|^2+2|a_2|^2)+|a_2|^4
-a_1\bar a_2^2a_3-\bar a_1a_2^2\bar a_3}.
\end{align*}

By the multi-point Schwarz-Pick lemma \eqref{eq:sp-n},
we have $|\gamma_n|=|g_n(0;0,\dots,0)|\le 1.$
Here, equality holds precisely if $g$ (equivalently $f$) is
a Blaschke product of degree $n.$
Schur \cite{Schur17} indeed showed that the sequence of inequalities
$|\gamma_n|\le 1$ characterizes the boundedness of an analytic function
$f$ by 1 in modulus.

Noting the relation $a_n=g^{(n)}(0)/n!=D_nf(z_0)/n!$ by Proposition
\ref{prop:rel}, we can rephrase the inequality $|\gamma_n|\le 1$
in terms of Peschl's invariant derivatives.
In particular, we obtain the following inequality
due to Yamashita as the case when $n=2.$

\begin{prop}[Yamashita $\text{\cite[Theorem 2]{Yam94P}}$]%
\label{lem:yamainequality}
Let $f\in\hol(\D).$
Then,
$$
|D_2f(z)|\le 2(1-|D_1f(z)|^2),\quad z\in\D.
$$
Equality holds for a point $z\in\D$ if and only if
$f$ is a Blaschke product of degree at most $2.$
\end{prop}

By the inequality $|\gamma_3|\le 1,$ we can similarly show the following.

\begin{thm}\label{thm:f3inequality}
Let $f\in\hol(\D).$
Then, for $ z\in\D,$
\begin{equation*}\label{eq:f3ineq}
\left| \frac{D_3f(z)}{6}\big(1-|D_1f(z)|^2\big)
+\overline{D_1f(z)}
\left(\frac{D_2f(z)}{2}\right)^2 \right|
+\left|\frac{D_2f(z)}2\right|^2
\le (1-|D_1f(z)|^2)^2,
\end{equation*}
where equality holds for a point $z\in\D$
if and only if $f$ is a Blaschke product of degree at most $3.$
\end{thm}

We have considered the simplest case when $z_0=z_1=\dots=z_j$ so far.
The second simplest case is perhaps when $z_0, z_1, \dots, z_j$ consist
of only two points.
We start with $f_2(z; z,z_0).$
The inequality $|f_2(z;z,z_0)|\le 1$ can be explicitly described
in the following result.

\begin{thm}[Generalized Dieudonn\'e's lemma]%\label{prop:d}
Let $f$ be an analytic function on $\D$
with $|f|<1$ and fix $z_0 \in \D.$
Then, for any point $z \in \D,$
\begin{equation*}\label{eq:f1ineq}
\left|f'(z)-\frac{f(z)-f(z_0)}{z-z_0} \cdot 
\frac{1-\overline{f(z_0)}f(z)}{1-|f(z_0)|^2}
\cdot \frac{1-|z_0|^2}{1-\overline{z_0}z} \right| 
\end{equation*} 
\begin{equation}\label{eq:f1ineq2}
\le \frac{1}{1-|z|^2}\left( \frac{|1-\overline{f(z_0)}{f(z)}|^2}{1-|f(z_0)|^2} 
\cdot
\left| \frac{z-z_0}{1-\overline{z_0}z} \right|
 - \frac{|f(z)-f(z_0)|^2}{1-|f(z_0)|^2} \cdot
\left| \frac{1-\overline{z_0}z}{z-z_0} \right| \right),
\end{equation} 
where equality holds if and only if $f$ is a Blaschke product of 
degree at most $2.$
\end{thm}

\begin{pf}
Let $g(z)=f_1(z;z_0).$
Then $f_2(z;z,z_0)=g_1(z;z).$
The inequality $|g_1(z;z)|=|D_1g(z)|\le 1$ is equivalent to
$$ |g'(z)| \le \frac{1-|g(z)|^2}{1-|z|^2}.$$
A straightforward calculation gives us the formula
$$
g'(z)=\frac{f'(z)}{[z,z_0]} \cdot 
\frac{1-|f(z_0)|^2}{(1-\overline{f(z_0)}f(z))^2}
-\frac{[f(z), f(z_0)]}{[z,z_0]^2} \cdot \frac{1-|z_0|^2}{(1-\overline{z_0}z)^2}.
$$
It takes a little rearrangements for the required inequality.
\end{pf}

Note that the inequality \eqref{eq:f1ineq2} is 
reduced to the original Dieudonn\'e's lemma when $z_0=f(z_0)=0:$
$$
\left|zf'(z)-f(z)\right|
\le \frac{|z|^2-|f(z)|^2}{1-|z|^2}.
$$
Conversely, through elementary computations, it can be seen that
the inequality \eqref{eq:f1ineq2} is obtained by applying
the original Dieudonn\'e's lemma to the function 
$h(\zeta)=[f([\zeta,-z_0]),f(z_0)]$
with the choice $\zeta=[z,z_0].$

It turns out that the inequalities $|f_2(z;z_0,z)|\le 1$ and
$|f_2(z_0;z,z)|\le 1$ are both equivalent to the inequality \eqref{eq:f1ineq2}.
Indeed, under the additional assumption that $z_0=f(z_0)=0,$
we have easily $|f_2(z;0,z)|=|f_2(0;z,z)|=\big|[D_1f(z),f(z)/z]/z\big|\le 1.$
If we set $w_0=[D_1f(z),f(z)/z]/z,$ we have $D_1f(z)=[zw_0,-f(z)/z]$
and the last inequality is equivalent to the
assertion $D_1f(z)\in A(\overline{\D}),$ where
$A(w)=[zw,-f(z)/z].$
Now use \eqref{eq:A} to see the equivalence with \eqref{eq:f1ineq2}.

We next consider the case when $j=3$ and $z_0, z_1, z_2, z_3$ consists of
two points.
By Lemma \ref{lem:inv}, we can assume that the two points are
$0$ and $z$ (which are interchangeable) and that $f(0)=0.$
The inequality $|f_3(z;0,0,0)|\le 1$ means that the third function
$f_3$ in the Schur algorithm has modulus at most 1.
On the other hand, the inequality $|f_3(z;z,z,0)|\le 1$
is rearranged to the following, which can be regarded as
Dieudonn\'e's lemma of the second order.

\begin{thm}
Let $f\in\hol(\D)\setminus\Aut(\D)$ with $f(0)=0.$
Then
\begin{align*}
&\left|
\frac{1}{2}z^2f''(z)-\frac{zf'(z)-f(z)}{1-|z|^2}
+\frac{\overline{f(z)}(zf'(z)-f(z))^2}{|z|^2-|f(z)|^2}
\right|
+\frac{|z||zf'(z)-f(z)|^2}{|z|^2-|f(z)|^2}
\\
&\le \frac{|z|(|z|^2-|f(z)|^2)}{(1-|z|^2)^2}.
\end{align*}
\end{thm}

\begin{pf}
Let $g(z)=f(z)/z.$
Then $g$ is a holomorphic self-map of $\D$
by assumption and $f_3(z;z,z,0)=g_2(z;z,z).$
The inequality $|f_3(z;z,z,0)|=|g_2(z;z,z)|\le 1$ is thus equivalent to
$\frac{1}{2}|D_2g(z)|+|D_1g(z)|^2\le 1$ 
(cf.~Proposition \ref{lem:yamainequality}).
The last inequality is indeed equivalent to the inequality in question.
\end{pf}

Finally, we consider the inequality $|f_3(z;z,0,0)| \le 1$ 
under the condition $f(0)=0.$    
Then we have the following inequality, which is another refinement of 
Dieudonn\'e's lemma involving the term $f'(0).$

\begin{thm}
Let $f\in\hol(\D)$ with $f(0)=0.$ Then
\begin{align*}
&\left|f'(z)(1-|f'(0)|^2)-\frac{2f(z)}{z}+  
\overline{f'(0)}\cdot \left(\frac{f(z)}{z}\right)^2 +f'(0) \right| \\
\le~ &\frac{1}{1-|z|^2}\left(\Big |z-\overline{f'(0)}{f(z)}\Big|^2-
\left| \frac{f(z)}{z}-f'(0) \right|^2  \right).
\end{align*} 
In particular, if $f'(0)=0$ in addition,
$$
\left|f'(z)-\frac{2f(z)}{z}\right| \le \frac{|z|^4-|f(z)|^2}{|z|^2(1-|z|^2)}.
$$
\end{thm}

\begin{pf}
Since $$f_3(z;z,0,0)=\frac{(1-|z|^2)f'_2(z;0,0)}{1- |f_2(z;0,0)|^2}$$ and
$$f_2(z;0,0)=\frac{1}{z} \cdot \frac{f(z)-zf'(0)}{z-\overline{f'(0)}f(z)},$$
we obtain the first inequality in the proposition 
by straightforward calculations.
\end{pf}

\def\cprime{$'$} \def\cprime{$'$} \def\cprime{$'$}
\providecommand{\bysame}{\leavevmode\hbox to3em{\hrulefill}\thinspace}
\providecommand{\MR}{\relax\ifhmode\unskip\space\fi MR }
\providecommand{\MRhref}[2]{%
  \href{http://www.ams.org/mathscinet-getitem?mr=#1}{#2}
}
\providecommand{\href}[2]{#2}

%\bibliography{papers}
\end{document}